\documentclass[12pt,english]{amsart}

\usepackage{amssymb, amsmath}
\usepackage{float,epsfig}
\usepackage{epsfig,amsfonts,amsbsy,amssymb,amsthm,here}
\usepackage[usenames,dvipsnames]{color}
\usepackage{graphics}
\usepackage{graphicx}
\usepackage{epsfig}
\usepackage{url}
\usepackage{arydshln,leftidx,mathtools}
\usepackage{rotating}
\usepackage{lscape}
\usepackage{setspace}
\usepackage{algpseudocode}
\usepackage{algorithm}
\usepackage[all]{xy}

\usepackage[paperwidth=220mm,paperheight=297mm,inner=3cm,outer=3cm,top=2.5cm,bottom=2.6cm]{geometry}

\newtheorem{theorem}{\bf Theorem}[section]
\newtheorem{proposition}[theorem]{\bf Proposition}
\newtheorem{definition}[theorem]{\bf Definition}
\newtheorem{corollary}[theorem]{\bf Corollary}
\newtheorem{conjecture}[theorem]{\bf Conjecture}
\newtheorem{example}[theorem]{\bf Example}

\newtheorem{remark}[theorem]{\bf Remark}
\newtheorem{lemma}[theorem]{\bf Lemma}

\newcommand{\ba}{\begin{array}}
\newcommand{\ea}{\end{array}}

\newcommand{\beano}{\begin{eqnarray*}}
\newcommand{\eeano}{\end{eqnarray*}}

\def \bmatrix#1{\left[ \begin{matrix} #1 \end{matrix} \right]}
\def\dmatrix#1{\left| \begin{matrix} #1 \end{matrix} \right|}
\newcommand \noin{\noindent}

\newcommand \R{{\mathbb R}}
\newcommand \C{{\mathbb C}}

\newcommand \V{{\mathcal V}}
\newcommand \I{{\mathcal I}}

\newcommand \tr{\mathrm{Tr}}

\newcommand \diag{\mathrm{Diag}}
\newcommand \eig{\mathrm{eig}}
\newcommand \pf{{\bf Proof: }}

\newcommand \x{{\mathbf{x}}}
\newcommand \y{{\mathbf{y}}}
\newcommand \z{{\mathbf{z}}}
\newcommand \0{{\mathbf{0}}}

\def \b{{\mathbf{b}}}
\def \tr{\mathrm{Tr}}

\def \u{{\mathbf{u}}}

\title{\textbf{Definite Determinantal Representations of  Multivariate  Polynomials}} 
\begin{document}
\author{Papri Dey}

\date{}
\maketitle
\begin{abstract}
In this paper, we consider the problem of  representing a multivariate polynomial  as the determinant of a definite (monic) symmetric/Hermitian linear matrix polynomial (LMP). Such a polynomial is known as determinantal polynomial. Determinantal polynomials can characterize the feasible sets of semidefinite programming problems that motivates us to deal with this problem. We introduce the notion of generalized mixed discriminant of matrices which translates the determinantal representation problem into computing a point of a real variety of a specified ideal. We develop an algorithm to determine such a determinantal representation of a bivariate polynomial of degree $d$. Then we propose a heuristic method to obtain a monic symmetric determinantal representation of a multivariate polynomial of degree $d$.

\end{abstract}
\noin \textbf{AMS Classification (2000)}. 14P99;  15A22; 15A39; 65F15; 65H04; 90C22. \\
\noin \textbf{Keywords.} Linear matrix Inequality, Linear Matrix Pencils, determinantal representation, RZ/Hyperbolic polynomials, semidefinite programming problems.

\section{Introduction}
A problem of characterizing multivariate polynomials which can be represented as the determinant of some monic (definite) symmetric/Hermitian linear matrix polynomial (LMP) is known as the determinantal representation problem in convex algebraic geometry. Although the constraint of monicness (definiteness) on LMP and coefficient matrices of a determinantal representation being symmetric/Hermitian make the problem more complicated, but it has generated a lot of interest due to its connection with \textit{semidefinite programming} (SDP) problems.

The technique of converting many types of optimization problems into SDP problems is successfully well-established and frequently
arise in control Theory, signal processing and many other areas in engineering. The
purpose of using this technique is that these converted problems can be solved by using
semidefinite programming algorithms. People use some popular tricks specific to certain
areas to convert it without having prior knowledge about whether this translation is possible theoretically.

Hence it is interesting to know which types of constraint sets are linear matrix inequality (LMI) representable sets.
To the best of author's knowledge this question was formally framed by Pablo Parrilo and
Bernd Sturmfels in \cite{Sturmfelsmin}. It is proved that
if the feasible set of an optimization problem is a definite LMI  representable set, the optimization problem can be transformed into a SDP problem \cite{Ramana1}, \cite{Helton}.

A set $S \subseteq \R^{n}$ is said to be \textit{LMI representable} if
\begin{equation} \label{lmiset}
S= \{ \x \in \R^{n} : A_{0} +x_{1}A_{1} +x_{2}A_{2} + \dots + x_{n}A_{n} \succeq 0 \}
\end{equation}
for some real symmetric matrices $A_{i}, i=0,\dots,n$ and $\x=(x_{1}, \dots, x_{n})^{T}$. If $A_{0}\succ 0$ (resp. $A_{0}=I$),
the set $S$ is called a \textit{definite} (resp. \textit{monic}) LMI representable set. 

It is proved in \cite{Helton} that if a polynomial $f(\x)$  is determinantal, then the algebraic interior associated with $f(\x)$ i.e., the closure of a (arcwise) connected component of $\{ \x \in \R^{n}:f(\x) > 0\}$  is a spectrahedron. Thus one of the successful techniques to deal with characterizing definite LMI representable sets is to characterize determinantal polynomials. 

In this paper, a polynomial $f(\x) \in \R[\x]$ is said to be determinantal if $f(\x)$ is the determinant of a definite (resp. monic) symmetric/Hermitian linear matrix polynomial (LMP); i.e.,
\begin{equation} \label{detreppoly}
f(\x) = \det(A_{0}+x_{1}A_{1} +x_{2}A_{2} + \dots + x_{n}A_{n}),
\end{equation}
where coefficient matrices  $A_{i}$ are symmetric/Hermitian of some order and the constant coefficient matrix $A_{0}$ is positive definite (resp. identity matrix).  The order of the coefficient matrices is called the size of determinantal representation and the size must  be  greater than or equal to $\deg(f)$.

A remarkable result due to Helton-Vinnikov \cite{Helton} proves that real zero (RZ) property of a polynomial is a necessary condition for the existence of monic symmetric/Hermitian determinantal representation. A polynomial $f(\x) \in \R[\x]$ with $f(0) \neq 0$ is said to be real zero (RZ) polynomial if its restriction along any line passing through origin has only real roots \cite{Helton}. The polynomial $f$  is called strictly RZ if all these roots are distinct, for all $\x \in \R^{n}, \x \neq 0$. 

Helton-Vinnikov have proved that RZ property is a necessary and sufficient condition for a bivariate polynomial to be a determinantal polynomial. The homogenized version of this result is known as Lax conjecture \cite{Pablo3}.  However, this result is no longer true for a  RZ polynomial in more than $2$ variables, i.e.,  it may not be a determinantal polynomial at all.  For example, dehomogenized polynomial of \textit{Vamos cube} $V_{8}$ is a RZ polynomial without a definite determinantal representation \cite{Branden}. This leads to the generalized Lax conjecture, for details see \cite{Vinnikov} \cite{Vinzant}, \cite{Pabloderivative}, \cite{Netzersmooth}.  Hyperbolic polynomials which play an important role in partial differential equations are the homegenized version of RZ polynomials \cite{Brahyper}.

The authors in \cite{Helton} have raised the attention towards computing such determinantal representations for any biariate polynomial and then this issue has been widely studied in literature, for example  one can see  \cite{Dixon}, \cite{Sturmfelsbivariate}, \cite{Henrion}, \cite{Vinnikov2}. To the best of authors' knowledge  nothing is known about multivariate determinantal polynomials except the fact that RZ property is a necessary condition for the existence of determinantal representation of a multivariate polynomial.

In this paper, we introduce the notion of generalized mixed discriminant of $k (\leq n)$ -tuple $n \times n$ matrices. Moreover, we prove that the coefficients of a multivariate determinantal polynomial can be uniquely determined by generalized mixed discriminant of the coefficient matrices of the corresponding determinantal representation, see Theorem \ref{themgmd}. 

We propose an algebraic method to compute the eigenvalues and diagonal entries of corresponding coefficient matrices for a determinantal multivariate polynomial. These results can be treated as necessary conditions for the existence of such a determinantal representation for a multivariate polynomial.

In fact, using the Theorem \ref{themgmd} we convert the determinantal representation problem of a bivariate polynomial into a problem of solving a system of polynomial equations, see Theorem \ref{themdet}. Moreover, we propose an algorithm to compute a monic symmetric determinantal representation of size $d$ for a bivariate polynomial of degree $d$ by finding a point in $\V_{\R}(\I)$. 

It is also shown that the ideal generated by these polynomials is generically a zero dimensional ideal and this method works efficiently for strictly RZ polynomials. The method has been exemplified for cubic and quartic bivariate polynomials symbolically and numerically. However, we study the equivalence classes of determinantal representations for bivariate polynomials, see Subsection \ref{equivclass}. 

We propose a heuristic method to determine a  monic symmetric determinantal representation of size $d$ for a multivariate polynomial of degree  $d$, see the Theorem \ref{themmulti}.

Note that a definite LMI representable set is always monic LMI representable \cite{Helton}.
So in this paper, we consider only problems dealing with monic symmetric/Hermitian determinantal representation of polynomials. 
We focus on the representations of size $d$ only. \textit{Monic symmetric determinantal representation} is abbreviated as  MSDR and  \textit{monic Hermitian determinantal representation} is abbreviated as MHDR in this paper.

\subsection*{Acknowledgements}
The author would like to thank her supervisor Prof. Harish K. Pillai for helpful discussions and suggestions on the subject of this paper. Much of the work on this paper has been supported by Council of Scientific and Industrial
Research (CSIR), India while the author was doing her doctoral work in IIT Bombay.
The author also gratefully acknowledges support through the Max Planck Institute for
Mathematics in the Sciences in Leipzig, Germany and the Institute for Computational
and Experimental Research in Mathematics, Brown University, USA.

\section{Generalized Mixed Discriminant of Matrices and Determinantal Polynomials}
 We introduce the concept of generalized mixed discriminant of $k (\leq n)$-tuple of $n \times n$  matrices (need not be distinct). Note that one can find the notion of mixed discriminant  of $n$-tuple of $n \times n$ distinct matrices \cite{Bapat} and  $k (< n)$ -tuple $n \times n$ distinct matrices in \cite{Srivastava}. Note that the  generalized mixed discriminant of matrices  is defined even if the matrices are not distinct. Then  by using the notion of  generalized mixed discriminant  we prove that the coefficients  of a determinantal multivariate polynomial can be  uniquely determined  in terms of the coefficient matrices of its determinantal representation.
\begin{definition} \label {generalizedmixdef}
Consider the $n \times n$ matrices $A^{(l)}=(a_{ij}^{(l)})$ for $l=1, \dots,n$. Pick any $k$-tuple matrix $(\underbrace{A^{(1)}, \dots , A^{(1)}}_{\mu_{1}}, \underbrace{A^{(2)}, \dots, A^{(2)}}_{\mu_{2}},\dots,\underbrace{A^{(n)},\dots,A^{(n)}}_{\mu_{n}}), \mu_{j}\in \{0,1,\dots,n\}, 1 \leq k \leq n$, and $\mu_{1}+\dots+\mu_{n}=k$. Then the generalized mixed discriminant (GMD) of the $k$ tuple of  $n \times n$ matrices is defined  as
\begin{equation} \nonumber
\widehat{D}(\underbrace{A^{(1)}, \dots , A^{(1)}}_{\mu_{1}}, \underbrace{A^{(2)}, \dots, A^{(2)}}_{\mu_{2}},\dots,\underbrace{A^{(n)},\dots,A^{(n)}}_{\mu_{n}}) =
\sum_{\alpha \in S[k]} \sum_{\sigma \in \widehat{S_{k}}(v)} 
   \dmatrix{a_{\alpha_{1}\alpha_{1}}^{(\sigma(1))} & \dots & a_{\alpha_{1}\alpha_{k}}^{(\sigma(1))} \\
     \vdots & & \\
     a_{\alpha_{k}\alpha_{1}}^{(\sigma(k))} & \dots & a_{\alpha_{k}\alpha_{k}}^{(\sigma(k))}}
\end{equation} where $S_{n}$ is the set of all permutations on $\{1, \dots,n\}$
  and $S[k]$ denotes the set of permutations of order $k$ which are chosen from the set $S_{n}$ such that $$\alpha = (\alpha_{1}, \dots, \alpha_{k}) \in S[k] \Rightarrow \alpha_{1} < \alpha_{2}< \dots < \alpha_{k},$$ $v=\{\underbrace{1, \dots, 1}_{\mu_{1}},\dots,\underbrace{n, \dots, n}_{\mu_{n}}\}$ and $\widehat{S_{k}}(v)$ is the set of all distinct permutations of $v$.
  \end{definition}

In order to derive the analytic expressions for the coefficients of a multivariate determinantal  polynomial in terms of
the coefficient matrices $A_{i}$s we need to prove the following results. 

\noin \textbf{Notation}: We follow the notation $|<\nabla m_{i1}(\x) \dots \nabla m_{in}(\x)>|$ to mean that the determinant of $M(\x)$ with $i$th row being replaced by $\nabla M_{i}(\x)$ where (the nabla symbol) $\nabla$ denotes the vector differential operator and $\x=(x_{1},\dots,x_{n})$.
\begin{lemma} \label{lemmamixcoeff}
Let $M(\x)=(m_{ij}(\x))$ be a $d \times d$ matrix with complex entries. Each entry of this matrix, denoted by   $m_{ij}(\x), i,j=1,\dots,d$  is a  linear polynomial in $x_{1},x_{2},\dots,x_{n}$ with complex coefficients.
Then
\begin{equation} \label{eqmix10}
\nabla |M(\x)| = \sum_{j=1}^{n}|< \nabla m_{j1}(\x) \dots \nabla m_{jn}(\x)>|.
\end{equation}
\end{lemma}

\pf  Without loss of generality assume that $M(\x):=A_{0}+x_{1}A_{1}+\dots+x_{n}A_{n}$. If $d=2$, then $M(\x):=\bmatrix{M_{1}\\M_{2}}$ and $\nabla |M(\x)|=|< \nabla  M_{1} >|+|< \nabla M_{2}>|$.
We prove this lemma by induction on the size $d$ of matrix $M(\x)$.

We assume that the equation (\ref{eqmix10})  is true  for all $M(\x) \in \C^{l \times l}[\x]$ . Now we have to show that this is true for $M(\x) \in \C^{(l+1) \times (l+1)}[\x]$.
Let $C_{1j}$ be the cofactor of $m_{1j}$. To derive the determinant of $M(\x)$ we consider the Laplace expansion of determinant along the Ist row of $M(\x)$. Then we have
\begin{align*}
&\nabla \dmatrix{M(\x)}=\sum_{j=1}^{l+1} \nabla (m_{1j} C_{1j})= \sum_{j=1}^{l+1} [C_{1j} \nabla m_{1j} +m_{1j} \nabla C_{1j}]  \\
& =\dmatrix{\nabla m_{11}(\x) \dots \nabla m_{1 (l+1)}(\x)}+m_{11} \sum_{i=2}^{l+1} \dmatrix{\nabla m_{i2}(\x) \dots \nabla m_{i (l+1)}(\x)} +\\
&\sum_{j=2}^{l+1} (-1)^{j-1}m_{1j} \sum_{k=2}^{l+1}\dmatrix{\nabla m_{k1}(\x) \dots \nabla m_{k(j-1)}(\x) \nabla m_{k (j+1)}(\x) \dots \nabla m_{k (l+1)}(\x)} 
\\
&=\dmatrix{\nabla m_{11}(\x) \dots \nabla m_{1 (l+1)}(\x)}+\sum_{j=2}^{l+1}\dmatrix{\nabla m_{j1}(\x) \dots
\nabla m_{j(l+1)}(\x)} \\
&=\sum_{j=1}^{l+1} |< \nabla m_{j1}(\x) \dots \nabla m_{jn}(\x)>|
\end{align*}
Thus it is true for $(l+1)$. So, by induction we can conclude that this result is true for any $d$. \qed

\begin{lemma} \label{lemmamixcoeffpar}
Let $M(\x)=A_{0}+x_{1}A_{1}+\dots+x_{n}A_{n}$ be a linear matrix polynomial of size $d$. For any $k_{j} \in \{0, \dots,d\}$ and $\sum_{j=1}^{n}k_{j}=l \leq d$
\begin{equation*}
\frac{\partial^{k_{1}+k_{2}+\dots+ k_{n}}}{\partial x_{1}^{k_{1}}x_{2}^{k_{2}}\dots x_{n}^{k_{n}}}\dmatrix{M(\x)}_{x_{1}=0,\dots,x_{n}=0}=(k_{1}! \dots k_{n}!) \widehat{D}(\underbrace{A_{0},\dots,A_{0}}_{d-l},\underbrace{A_{1},\dots,A_{1}}_{k_{1}} \dots \underbrace{A_{n},\dots,A_{n}}_{k_{n}})
\end{equation*}
\end{lemma}
\pf  By Lemma \ref{lemmamixcoeff} it is clear that  L.H.S  is equal to the sum of determinants of some matrices. These matrices are constructed as follows. In the Ist  partial derivative of $M$ with respect to $x_{j}$, new matrices are constructed such that  only one row of $M$ is replaced by the  corresponding row of the  matrix  $A_{j}$. Using the same logic we claim that in the $k_{1}k_{2}\dots k_{n}$-th derivative of $M$ , the $k_{1}$ rows of  $M$ will be replaced by the $k_{1}$ rows of $A_{1}$, $k_{2}$ rows of $M$ will be replaced by $A_{2}$ and $k_{n}$ rows will be replaced by $k_{n}$  rows of $A_{n}$. 
If $k_{1}+\dots+k_{n}=l <d$, then  the remaining $d-l$ rows of $M$ will be replaced by the corresponding rows of matrix $A_{0}$  at  $x_{1}=\dots=x_{n}=0$.  From the Definition \ref {generalizedmixdef}   it is clear that  we need to construct such matrices in order to calculate  generalized mixed discriminant of matrices. So, the partial derivatives of the determinant of any multivariate linear matrix  polynomial  can be determined by calculating the generalized mixed discriminant of certain matrices. 

Now we need to prove that the number of determinants in both sides of the equality is same.  As there are $d$ rows in the linear matrix polynomial $M$, so there are $d$ determinants in  $\frac{\partial}{\partial x_{j }}\dmatrix{M(x_{1},\dots,x_{n})}$. If we differentiate it one more time, each determinant provides $(d-1)$ nontrivial determinants. So, the number of nontrivial determinants in  L.H.S  is $d(d-1) \dots (d-k_{1}+1) \dots (d-k_{1}-\dots-k_{n}+1)$. Due to Leibnitz product  rule for differentiation there is a functional equality or parity among  determinants in the expansion of the L.H.S term. Thus there are repeated determinants in the expansion  and each distinct determinant must be repeated the same number of times; i.e.; $k_{1}!\dots k_{n}!$. On the other hand, in order to calculate the generalized mixed discriminant of matrices we follow the lexicographic order;  i.e., ($i_{j} \ < i_{k}$ if $ j < k$ ) to choose $k_{1}, \dots,k_{n}$  rows out of $d$ rows. So, there are $d \choose k_{1}$$d-k_{1} \choose  k_{2}$$ \dots$$ d-k_{1}-k_{2}-\dots-k_{n-1} \choose k_{n}$ determinants in $ \widehat{D}(\underbrace{A_{0},\dots,A_{0}}_{d-l},\underbrace{A_{1},\dots,A_{1}}_{k_{1}} \dots \underbrace{A_{n},\dots,A_{n}}_{k_{n}})$ . It satisfies the following identity
\begin{equation*}
d(d-1)\dots (d-k_{1}-k_{2}-\dots-k_{n}+1)=k_{1}!\dots k_{n}! \mbox{$d \choose k_{1}$}\mbox{$d-k_{1} \choose  k_{2}$}\mbox{$ \dots$}\mbox{$ d-k_{1}-k_{2}-\dots-k_{n-1} \choose k_{n}$}.
\end{equation*}
Hence the desired result follows. \qed

\begin{remark} {\rm
Similar kinds of results exist in literature  \cite{Alexandria},  \cite{Bhatia}, \cite{Schneidermixed}.
}
\end{remark}
\begin{theorem} (Generalized Mixed Discriminant Theorem) \label{themgmd}
The coefficients of a multivariate determinantal olynomial $f(\x) \in \R[\x]$ of degree $d$ are uniquely determined by the generalized mixed discriminants of the coefficient matrices $A_{i}$ as follows. If the degree of a monomial $x_{1}^{k_{1}}x_{2}^{k_{2}} \dots x_{n}^{k_{n}}$ is $k (k_{1}+k_{2}+\dots+k_{n}=k) \leq d$, then the coefficient $f_{k_{1}\dots k_{n}}$ of
 ($x_{1}^{k_{1}} \dots x_{n}^{k_{n}})$ is given by $$\widehat{D}(\underbrace{A_{1}, \dots , A_{1}}_{k_{1}}, \underbrace{A_{2}, \dots, A_{2}}_{k_{2}},\dots,\underbrace{A_{n},\dots,A_{n}}_{k_{n}}).$$
In particular,
\begin{enumerate}
\item $\tr A_{i} =$ the coefficient of $x_{i}$ for all $i=1,\dots,n$.
\item $\det A_{i} =$ the coefficient of $x_{i}^{d}$ for all $i=1,\dots,n$ where $d$, the degree of the polynomial is equal to the size $d$ of the coefficient matrix.
\end{enumerate}
\end{theorem}
\pf  As $f(\x)$ is a determinantal polynomial, so  $f(\x)=\det(I+x_{1}A_{1}+\dots+x_{n}A_{n})$ for some symmetric matrices $A_{i},i=1(1)n$.
Using multivariate Taylor series coefficient formula the coefficients of $\det(I+x_{1}A_{1}+\dots+x_{n}A_{n})$ can be determined by the given formula \cite{Apostol}.
\begin{equation*} \label{eqmix12}
f_{k_{1}k_{2} \dots k_{n}}=\frac{1}{k_{1}!k_{2}! \dots k_{n}!}\frac{\partial^{k_{1}+k_{2}+\dots+ k_{n}}}{\partial x_{1}^{k_{1}}x_{2}^{k_{2}}\dots x_{n}^{k_{n}}} |I+x_{1}A_{1}+\dots+x_{n}A_{n}|_{x_{1}=\dots=x_{n}=0}, k_{j} \in \{0,\dots, d\} 
\end{equation*}
Therefore by Lemma \ref{lemmamixcoeffpar} we  conclude that
\begin{equation*} \label{eqmix12}
f_{k_{1}k_{2} \dots k_{n}}=\widehat{D}(\underbrace{A_{1}, \dots , A_{1}}_{i_{1}}, \underbrace{A_{2}, \dots, A_{2}}_{i_{2}},\dots,\underbrace{A_{n},\dots,A_{n}}_{i_{n}})
\end{equation*}
\qed
\section{Determinantal Polynomials}
At first we discuss some facts about determinantal multivariate polynomials. 
Since the coefficient matrices $A_{i}$s are Hermitian (symmetric), therefore by the
spectral theorem of a Hermitian (symmetric) matrix there exist a unitary (orthogonal) matrix $U_{i}$
such that $A_{i}=U_{i}^{*}D_{i}U_{i}$ for all $i=1, \dots,n$ where $D_{i}$ is a diagonal matrix  whose diagonal entries are the eigenvalues of $A_{i}$. So, one can always find a suitable unitary (orthogonal) matrix $U$ such that one of the coefficient matrices becomes diagonal. Without loss of generality, it is enough to consider coefficient matrix $A_{1}$ associated to $x_{1}$ as a diagonal matrix $D_{1}$ and obtain an MSDR (MHDR) of the following form
\begin{equation} \label{eqdet}
f(\x)=\det(I+x_{1}D_{1}+x_{2}V_{12}D_{2}V_{12}^{\ast}+\dots+x_{n}V_{1n}D_{n}V_{1n}^{\ast})=\det(I+x_{1}D_{1}+x_{2}A_{12}\dots+x_{n}A_{1n})
\end{equation}
 where $V_{ij}, i \neq j$ is the transition matrix from $D_{i}$ to $A_{ij}:= V_{ij} D_{j}V_{ij}^{\ast}$
 and similarly $V_{ij}^{\ast}=V_{ji}, i \neq j$ is the transition matrix from $D_{j}$ to $ A_{ji}:=V_{ij}^{T} D_{i}V_{ij}=V_{ji}D_{i}V_{ij}$.

\subsection{Eigenvalues of Coefficient Matrices}
Observe that the eigenvalues of the coefficient matrices $A_{1i}$ are nothing but the entries of the diagonal matrices $D_{i}$ for all $i=2,\dots,n$. We explain a technique to determine these diagonal matrices. 

We take restrictions of the given multivariate polynomial $f(\x)$ along each  $x_{i}, i=1, \dots,n$ that means we  restrict the polynomial along one variable at a time by making the rest of the variables zero and generate $n$ univariate polynomials $f_{x_{i}}=f(0,\dots,x_{i},\dots,0)$. 

It is known that if a multivariate polynomial $f(\x)$ admits an MSDR (MHDR), it is a RZ polynomial. By recalling the definition of RZ polynomial, we know that for any $\x \in \R^{n}$, RZ polynomial $f(\x)$ when restricted along any line passing through origin, has only real zeros. So when a RZ polynomial $f(\x)$ restricted along $x_{i}, i=1,\dots,n$, each of them has only real zeros, i.e., all univariate polynomials $f_{x_{i}}$ in $x_{i}$ have only real zeros.

As a consequence of this result we have a necessary condition for the existence of an MSDR of size equal to the degree of the polynomial for a multivariate polynomial of any degree.
\begin{lemma}\label{existence of diagonal 2}
If a multivariate polynomial $f(\x) \in \R[\x]$ of degree $d$ has an MSDR (MHDR) of size $d$, then all the roots of $f_{x_{i}}$ are real for all $i=1,\dots,n$.
\end{lemma}

More interestingly, the entries of the diagonal matrices $D_{i}$ can be obtained from the roots of $f_{x_{i}}$ for all $i=1, \dots, n$ by the following Lemmas.
\begin{lemma} \label{reverse}
For a univariate polynomial $f(x)=\det(xI+A)$ of degree $d$, the reversed linear polynomial $\widehat{f}(\x):=x^{d}f(1/x) =\det(I+xA)$ has the same coefficients in the reverse order.
\end{lemma}
\pf Say, $f(x)=\sum_{i=0}^{d}a_{d-i}x^{i}, a_{0} =1$, then  $\widehat{f}:=x^{d} f(1/x)=\det(I+xA)= \sum_{i=0}^{d}a_{i}x^{i},a_{0}=1$.
Since there is a one to one correspondence between the roots of $f$ and $\widehat{f}$ which is given by $x \rightarrow 1/x$, therefore the roots of $\widehat{f}$  are the reciprocals of the nonzero finite eigenvalues of $A$. If $A$ is a singular matrix i.e., it has zero eigenvalues, then the degree of $\widehat{f}$ is dropped by the number of multiplicities of the zero eigenvalues of $A$, but the coefficients of these two polynomials are same in the reverse order. \hfill{$\square$}

Note that
\begin{equation*}
f(x)=\det(xI+A)=x^{d}+E_{1}(A) x^{d-1} + \dots + E_{d}(A)
\end{equation*}
where $E_{k}(A)$ is the sum of $k \times k$ principal minors of $A, k=1, \dots, d$. If $\lambda_{1}, \dots, \lambda_{d}$ are the eigenvalues of $A \in M^{d \times d}$, then the sum of $k \times k$ principal minors of $A$ is the $k$-th elementary symmetric function of the eigenvalues of $A$ i.e., $S_{k}(\lambda_{1}, \dots, \lambda_{d}) = E_{k}(A)$ \cite{Johnson}. The $k$-th elementary symmetric function of $d$ numbers $\lambda_{1}, \dots, \lambda_{d}, k \leq d$ is defined as $S_{k}(\lambda_{1}, \dots, \lambda_{d}) = \sum_{1 \leq t_{1} < \dots < t_{k}=d } \prod_{j=1}^{k} \lambda_{t_{j}}$. In particular,
$\tr(A)= \sum_{i=1}^{d} \lambda_{i}, \det(A)= \prod_{i=1}^{d} \lambda_{i}$. Thus
\begin{equation*}
\widehat{f}(x)=\det(I+xA)=E_{d}(A)x^{d}+E_{d-1}(A)x^{d-1}+\dots +1
\end{equation*}
\begin{lemma} \label{lemmartdiag}
The non-zero eigenvalues of coefficient matrices of a determinantal polynomial defined in equation (\ref{eqdet}) are the negative reciprocal of the roots of univariate polynomials $f_{x_{i}}:=f(0,\dots, x_{i},\dots,0)$ for all $i=1,, \dots,n$.
\end{lemma}
\pf By Lemma \ref{reverse} a univariate polynomial $f(t)$ has only real zeros and $f(0) \neq 0$ if and only if the reversed polynomial $t^{d}f(1/t)$ has only real zeros. So, a polynomial $f(\x)$ is a RZ polynomial if and only if for any fixed real vector $\x \in \R^{n}$ the univariate polynomial $\widehat{f}_{\x}(t):=t^{d}f(\x / t)$ in $t$ has only real zeros. Thus, if a polynomial $f(\x) \in \R[\x]$ is a RZ polynomial, then the associated univariate polynomial of the homogenized polynomial $\widehat{f_{x_{i}}}(t):=t^{d}f_{x_{i}}(x_{i}/t)$
has only real zeros at point $\x=(0,\dots,x_{i},\dots,0) \in \R^{n}$. As $\det(tI+D_{i})=\widehat{f_{e_{i}}}(t)$, here map is:$t \mapsto -t$, therefore the roots of the reversed polynomials $\widehat{f_{e_{i}}}(t):=t^{d}f_{e_{i}}(e_{i}/t)$ at $\x=e_{i}$ are precisely the negative of the diagonal elements of $D_{i}$ for all $i=1, \dots, n$. As polynomials $f_{x_{i}}$ can be viewed as the reversed polynomials of $\widehat{f_{e_{i}}}(t)$ in one parameter $t$ instead of $x_{i}$, so the entries of diagonal matrices $D_{i}$ are the negative reciprocal of the roots of univariate polynomials $f_{x_{i}}$ for all $i=1, \dots,n$. As the coefficient matrices are either symmetric or Hermitian, so by the spectral theorem of symmetric or Hermitian matrices, $\eig(A_{1i})=\diag(D_{i})$, for all $i=2,\dots, n$. \qed


\begin{remark}\rm{
There are many ways to calculate the roots of a univariate polynomial. One of the popular methods is based on using the companion matrix associated to that polynomial. The eigenvalues of the companion matrix $C_{\widehat{f_{e_{i}}}(t)}$ associated with polynomial $\widehat{f_{e_{i}}}(t)$ are
the roots of the polynomial $\widehat{f_{e_{i}}}(t)$ since $\det(tI-C_{\widehat{f_{e_{i}}}}(t))=\widehat{f_{e_{i}}}(t)$.}
\end{remark}
Note that diagonal entries of $D_{j}$ are actually the eigenvalues of $A_{1j}$ defined in equation (\ref{eqdet}).
\subsection*{The diagonal entries of coefficient matrices}
Find the diagonal entries of coefficient matrix $A_{1i}$ for all $2 \leq i \leq n$.
Say $D_{1}=\diag(r_{1},r_{2},\dots,r_{d})$   and $D_{2}=\diag(s_{1},s_{2},\dots,s_{d})$. 

Using the  generalized mixed discriminant of coefficient matrices, by the Theorem \ref{themgmd} we derive the analytic expressions for each coefficient of  $f(\x)$ in terms of  the entries of $D_{1}$ and $A_{1i}$. We obtain analytic expressions for the vector coefficient  of mixed monomials $x_{1}^{\alpha_{1}}x_{i}^{\alpha_{2}},$
where $1 \leq \alpha_{1} \leq d-1$ and $1 \leq \alpha_{2} \leq d-1$ by the Theorem  \ref{themgmd}.

The diagonal entries of a coefficient matrix $A_{1i}$  can be determined  by solving a system of linear equations of the form $G \y_{i} =\z_{i}$  for all $i=2,\dots,n$  where $\y_{i}$  denotes the vector consisting of the diagonal entries of $A_{1i}$, i.e., $\y_{i}:=\diag(A_{1i})$, and
 \begin{equation} \label{eqdiag}
 \z_{i}=\bmatrix{\mbox{coeff of} \ x_{i} \\ \mbox{coeff of} \ x_{1}x_{i}\\ \vdots  \\ \mbox{coeff of}\  x_{1}^{d-1}x_{i}}, G=\bmatrix{1 & 1 & \dots & 1\\ \sum_{i=2 }^{d}r_{i}& \sum_{i=1, i \neq 2}^{d}r_{i}& \dots &  \sum_{i=1}^{d-1}r_{i}\\ \sum_{i_{k},i_{l}\neq 1,i_{k} <i_{ l}}r_{i_{k}}r_{i_{l}} & \dots & \dots & \sum_{i_{k},i_{l}\neq d,i_{k} <i_{ l}}r_{i_{k}}r_{i_{l}}  \\ \vdots & \vdots & \vdots & \vdots \\ r_{2}r_{3}\dots r_{d} & r_{1}r_{3}\dots r_{d} & \dots & r_{1}\dots r_{d-1}}.
\end{equation}

As $\z_{i}$ is the vector of coefficients of monomials of the form $x_{1}^{\alpha_{1}}x_{i}, 0 \leq \alpha_{1} \leq d-1$, so 
the relations are linear in terms of the entries of diagonal entries of $A_{1i}$. Thus by the Theorem \ref{themgmd} $\diag(A_{1i})$ can be found by solving a system of linear equations. 

Note that the number of mixed monomials which are of the form $x_{1}^{\alpha_{1}}x_{i}$ with
$0 \leq \alpha_{1} \leq d-1$ is $d$  and the diagonal entries of coefficient matrix $A_{1i}$ is also $d$. So, if the polynomial $f(\x)$ is determinantal, i.e., it satisfies the equation (\ref{eqdet}), one can always determine the diagonal entries of coefficient matrix $A_{1i}$ by solving a system of $d$ linear equations in $d$ unknowns. However, they may not be unique.
\begin{lemma} \label{lemmauniquediag}
A matrix $G$ defined in equation (\ref{eqdiag})
is invertible if and only if diagonal entries of $D_{1}$ defined in equation (\ref{eqdet}) are all distinct.
\end{lemma}
\pf Note that $D_{1}=\diag(r_{1}, \dots, r_{d})$. The result follows as $\det(G)=\prod_{i,j=1,i <j}^{d}(r_{i}-r_{j})$. \qed

As a result, we have
\begin{corollary}
The diagonal entries of $A_{1i}$ are uniquely determined up to ordering if all the diagonal entries of coefficient  matrix $D_{1}$ are distinct.
\end{corollary}
\begin{remark}\rm{
 Without loss of generality one can make any of the coefficient matrices as diagonal matrix, so by Lemma \ref{lemmauniquediag} the diagonal entries of coefficient matrices of a determinantal multivariate  polynomial can be uniquely determined if and only if one of the coefficient matrices have all distinct eigenvalues.
 }
\end{remark}
\section{Bivariate Polynomials}
In this section, we propose  a method to determine an MSDR (MHDR) of size $d$ for a bivariate polynomial of degree $d$ by solving a system of polynomial equations. We would like to mention that though the  authors in \cite{Sturmfelsbivariate} have proposed a method to determine an  MSDR  by solving polynomial equations, but note  that they have not  explained any method to obtain the analytic expressions of the coefficients of a determinantal bivariate polynomial in terms of corresponding coefficient matrices. 

In this paper, we propose an explicit method to express the coeffcients of a determinantal  multivariate polynomial  in terms of coefficient matrices  of the corresponding determinantal representation of that polynomial and convert the determinantal representations problem into solving a system of polynomial equations.

By equation (\ref{eqdet}) a determinantal bivariate polynomial $f(\x) \in \R[\x]$ can be written as
\begin{equation} \label{eqbivdet}
 f(\x)=\det(I+x_{1}D_{1}+x_{2}V_{12}D_{2}V_{12}^{T})=det(I+x_{1}D_{1}+x_{2}A_{12})
\end{equation}

By Lemma \ref{lemmartdiag} one can find $\eig(D_{1})$ and $\eig(A_{12})$. The diagonal entries of matrix $A_{12}$ in equation (\ref{eqbivdet}) can be determined  by solving a system of linear equations of the form $G \y_{2} =\z_{2}$   where $\y_{2}=\diag(A_{12})$, and 
 $\z_{2}$, $G$ are defined in equation (\ref{eqdiag}).
 
As the diagonal entries of $A_{12}$ are evaluated, so the number of unknown entries of $A_{12}$ which are the off-diagonal entries is $(d+1) \choose 2$ $-d=$ $d \choose 2$


\subsection*{Computation of MSDR (MHDRs)}
Consider the expressions associated with the coefficients of monomials
$x_{1}^{\alpha_{1}}x_{2}^{\alpha_{2}}, 0 \leq \alpha_{1} \leq d-2, 2 \leq \alpha_{2} \leq d$ by the Theorem \ref{themgmd}. 
The number of monomials of the the form  $x_{1}^{\alpha_{1}}x_{2}^{\alpha_{2}},$
where $0 \leq \alpha_{1} \leq d-2$ and $2 \leq \alpha_{2} \leq d$ is $(d-1+d-2+d-3+\dots+2+1)=$ $d \choose 2$.
Then by comparing the expressions associated with the coefficients of the remaining monomials $x_{1}^{\alpha_{1}}x_{2}^{\alpha_{2}},$
where $0 \leq \alpha_{1} \leq d-2$ and $2 \leq \alpha_{2} \leq d$, we obtain generically a zero dimensional ideal
$\I$ generated by  $d \choose 2$ polynomial equations in $d \choose 2$ parameters.

Note that if any two diagonal entries of $D_{1}$ or $D_{2}$ are equal, ideal $\I$ may or may not be a zero dimensional ideal. As the diagonal entries of $D_{1}$ and $A_{12}$ are known and diagonal entries of coefficient matrix $D_{1}$ are generically distinct, so we obtain generically a zero dimensional ideal $\I$ generated by $d \choose 2$ polynomials in $d \choose 2$ parameters.

Find an element of the  real variety $\V_{\R}(\I)$. This can be obtained by using available softwares like
\textbf{Singular}, \textbf{Bertini}, \textbf{Maple} and \textbf{Sage}. For instance,  using Sage find Groebner basis of the ideal $\I$ and then find real roots of that Groebner basis. Check whether there exists at least one real root, otherwise exit-no MSDR (MHDR) of size $d$ is possible.  Find the off diagonal entries of coefficient matrix $A_{12}$. Therefore, 
we conclude the following theorem.

\begin{theorem} \label{themdet}
  A bivariate polynomial of degree $d$ is determinantal i.e., $f(\x)=\det(I+x_{1}D_{1}+x_{2}A_{12})$ if and only if $\V_{\R}(\I)$ is non-empty where $\I$ is generated by $d \choose 2$ polynomials with real coefficients obtained by Theorem \ref{themgmd} associated with the coefficients of monomials $x_{1}^{\alpha_{1}}x_{1}^{\alpha_{2}}, 0 \leq \alpha_{1} \leq d-1,2 \leq \alpha_{2} \leq d$ in $d \choose 2$ unknowns which are the off-diagonal entries of coefficient matrix $A_{12}$.
\end{theorem}

Thus we propose the following  algorithm to determine an MSDR of size $d$ for bivariate polynomial by solving a system of polynomial equations.
\begin{algorithm}\label{algobiv}
  \caption{Algorithm to Determine an MSDR of size $d$ for Bivariate Polynomial of degree $d$}
\begin{algorithmic}
\State Input: Bivariate polynomial $f(\x)$ of degree $d$.
\State Output: Coefficient matrices $D_{1}, A_{12}$ such that
\begin{equation*}
f(\x)=\det(I+x_{1}D_{1}+x_{2}A_{12}).
\end{equation*}
\\\hrulefill
\begin{enumerate}
\item Construct the univariate polynomials $f_{x_{i}}, i=1,2$.
\item Determine the eigenvalues of coefficient matrices $D_{1},A_{12}$ by Lemma \ref{lemmartdiag}.
\item Check that the  eigenvalues of coefficient matrices $D_{1},A_{12}$ are real. If not, exit-no MSDR of size $d$ possible.
\item Find the diagonal entries of the matrix $A_{12}$ by solving a system of linear equations of the form $G \y_{2}=\z_{2}$ defined in equation (\ref{eqdiag}).
\item Comparing the expressions associated with the coefficients of monomials $x_{1}^{\alpha_{1}}x_{2}^{\alpha_{2}},$
where $0 \leq \alpha_{1} \leq d-2$ and $2 \leq \alpha_{2} \leq d$, by the Theorem \ref{themgmd} obtain a generically zero dimensional ideal
$\I$ generated by $d \choose 2$ polynomial equations in $d \choose 2$ parameters.
\item Find an element of real variety $V_{\R}(\I)$. 
\item Check whether there exists at least one  REAL root, otherwise exit-no MSDR of size $d$ possible.
\item Construct $D_{1}$ and $A_{12}$.
\end{enumerate}
\end{algorithmic}
\end{algorithm}
\subsection{Method of Solving Polynomial Equations}
In this subsection, we study the method of solving polynomial equations for cubic and quartic bivariate case. For cubic bivariate case we explain the method symbolically as well as numerically. Consider the cubic bivariate polynomial 
 \begin{equation} \label{cubicbiveq1}
 f(x_{1},x_{2})=
 f_{30}x_{1}^{3}+f_{03}x_{2}^{3}+f_{21}x_{1}^{2}x_{2}+f_{12}x_{1}x_{2}^{2}+f_{20}x_{1}^{2}+f_{02}x_{2}^{2}+f_{11}x_{1}x_{2}+f_{10}x_{1}+f_{01}x_{2}+1.
 \end{equation}

Say the coefficient matrix $A_{12}=\bmatrix{a & d & e \\d & b & f \\e & f & c}$. If MSDR of size $3$ exists for a given cubic bivariate polynomial, by equation (\ref{eqbivdet}) it is of the form
\begin{equation} \label{MSDReq2}
f(\x)=:\det(I_{3}+x_{1}\bmatrix{d_{1} & 0 & 0\\0 & d_{2} & 0 \\ 0 & 0 & d_{3}}+x_{2}\bmatrix{a & d & e \\d & b & f \\e & f & c})
\end{equation}
Determine the diagonal entries of $D_{1}$ by Lemma \ref{lemmartdiag} and check whether they are real. Observe that $D_{1}$ is unique up to ordering $d_{1} > d_{2} > d_{3}$.

By solving a system of linear equations of the form $G \y =\z$ where
\begin{equation*}
G=\bmatrix{1 & 1 & 1\\(d_{2}+d_{3})&(d_{1}+d_{3})&(d_{1}+d_{2})\\d_{2}d_{3}&d_{1}d_{3}&d_{1}d_{2}}, \y=\bmatrix{a\\b\\c},\z=\bmatrix{f_{01}\\ f_{11} \\ f_{21}},
\end{equation*}
we obtain $a,b$ and $c$ provided the system is consistent. If the system is inconsistent, then no MSDR of size $3$ is possible.

Two of these three equations are quadratic and other one is cubic.

By the Theorem \ref{themdet} the ideal $\I$ is generated by polynomial equations associated with monomials $x_{2}^{2},x_{1}x_{2}^{2},x_{2}^{3}$.  We derive the following polynomial equations with real coefficients (as $d_{1},d_{2},d_{3},a,b,c$ are real) in three unknowns $d,e,f$ by the Theorem \ref{themgmd}. This can also be obtained by expanding and comparing the coefficients of equation (\ref{MSDReq2}).
\begin{align} \label{coefftensoreq3}
&f_{02}=bc+ac+ab-f^{2}-d^{2}-e^{2} \nonumber \\
&f_{12}=d_{1}bc+abd_{3}+ad_{2}c-d_{1}f^{2}-d_{3}d^{2}-d_{2}e^{2} \nonumber \\
&f_{03}=abc+2def-af^{2}-cd^{2}-be^{2} \nonumber 
\end{align}
We obtain the ideal
\begin{align*}
\I=&\{d^{2}+e^{2}+f^{2}-(bc+ac+ab)+ f_{02}, d_{3}d^{2}+d_{2}e^{2}+d_{1}f^{2}-(d_{1}bc+d_{2}ac+d_{3}ab)+ f_{12}\\&, cd^{2}+be^{2}+af^{2}-(abc+2def)+ \ f_{03}\}
 \end{align*}
which is generically a zero dimensional ideal generated by $3$ polynomial equations in $3$ unknowns. Calculate Groebner basis of the ideal $\I$ (preferably in SAGE). Use Solve command to find  real roots of the Groebner basis to find $d,e,f$. Determine the off-diagonal entries of $A_{12}$. Check whether there exists at least one  REAL root, otherwise no MSDR of size $3$ is possible. Thus we can construct $D_{1}$ and $A_{12}$.

For cubic bivariate case we propose an alternative method here to compute such a determinantal representation instead of finding $\V_{\R}(\I)$. 

Observe that the equations due to coefficients $x_{2}^{2},x_{1}x_{2}^{2}$ are linear in $d^{2},e^{2},$and $f^{2}$. So, we can determine the set $\mathcal{K}$ of solutions of this under determined system of linear equations $ \bmatrix{1 & 1 & 1\\ d_{3} &d_{2}& d_{1}} \bmatrix{d^{2} \\e^{2}\\f^{2}}=\bmatrix{(bc+ac+ab)-f_{02}\\ (d_{1}bc+d_{2}ac+d_{3}ab)-f_{12}}$ i.e.,
\begin{equation*}
\mathcal{K}:=\{\bmatrix{(d^{2})^{\ast}\\(e^{2})^{\ast}\\(f^{2})^{\ast}}+k \gamma :\gamma \in \ker(\bmatrix{1 & 1 & 1\\ d_{3} &d_{2}& d_{1}}) \}
\end{equation*}
where $\bmatrix{(d^{2})^{\ast}\\(e^{2})^{\ast}\\(f^{2})^{\ast}}$ is a solution of the under determined system. Using the equation of $x_{2}^{3}$ one can derive a cubic equation in $k$ and by solving that equation one can obtain MSDRs if it exists.

As we know from linear algebra that if $\z \in$ col$(G),$ column space of $G$, the system is consistent. Here we have to study three cases separately.
\begin{itemize}
\item Diagonal matrices $D_{1},D_{2}$ are simple i.e., all three entries of each of diagonal matrices $D_{i}, i=1,2$ are distinct:
Note that $G$ is invertible in this case. So, the system has a unique solution. So, the diagonal entries of coefficient matrix $A_{12}:=V_{12}D_{2}V_{12}^{T}$ are uniquely determined in this case. Due to symmetry in coefficients of a polynomial we could make the coefficient matrix associated to $x_{2}$ a diagonal matrix and similarly we had to determine the diagonal entries of coefficient matrix $V_{12}^{T}D_{1}V_{12}$. As it is assumed that both of the diagonal matrices are simple, so the diagonal entries of coefficient matrix $V_{12}^{T}D_{1}V_{12}$ are uniquely determined by the same logic.
\item At least one of the two diagonal matrices is having the following property: (All three entries are equal)

As without loss of generality (wlog) we can make any one of two coefficient matrices a diagonal matrix, so we can choose any of these two coefficient matrices have the mentioned property. Wlog say diagonal matrix $D_{1}$ has this property, so $D_{1}=\lambda I_{3}$, identity matrix of order $3$ and $\lambda$ is a non zero scalar (otherwise diagonal matrix would be a zero matrix). Observe that matrix $G$ is not invertible in this case. If solution exists, there are infinitely many solutions in this case. This result reflects the fact that there are infinitely many (orthogonally equivalent) symmetric representations as
\begin{equation*}
f(x_{1},x_{2})=\det(I+x_{1} \lambda I +x_{2} D_{2})=\det(I+x_{1} \lambda I +x_{2}VD_{2}V^{T})
\end{equation*}
for any orthogonal matrix $V$ of order $3$.
\item At least one of the two diagonal matrices is having the following property: (Two of three entries are equal)

In this case the ideal $\I$ may or may not be a zero dimensional ideal. If it is not a zero dimensional ideal \cite{Cox1}, then the alternative method works well as opposed to finding $\V_{\R}(\I)$.
\end{itemize}

\begin{example}\label{excubicbivtensor}{\rm
Consider the cubic bivariate polynomial $f(\x)=6x_{1}^3 + 36x_{1}^2x_{2} + 11x_{1}^2 + 66x_{1}x_{2}^2 + 42x_{1}x_{2} + 6x_{1} + 36x_{2}^3 + 36x_{2}^2 + 11x_{2} + 1$. So, the univariate polynomial $f_{x_{1}}:=6x_{1}^{3}+11x_{1}^{2}+6x_{1}+1$. By calculating the roots of $f_{x_{1}}$, determine $D_{1}=\diag(3,2,1)$. Derive a system of three linear equations in three unknowns $a,b$ and $c$ due to the relations of coefficients of $x_{2},x_{1}x_{2},x_{1}^{2}x_{2}$. By solving this system of linear equations we obtain the diagonal entries of $A_{12}$. Here we have $\bmatrix{a\\b\\c}=\bmatrix{9/2\\4\\5/2}$ which are diagonal entries of $A_{12}$. Now we need to determine the off-diagonal entries of $A_{12}$. As we know the values of $d_{1},d_{2},d_{3}$ and $a,b,c$, so we can simplify the equations in (\ref{coefftensoreq3}) and write those analytic expressions in terms of off-diagonal entries $d,e,f$ and calculate the constant terms from the polynomial obtained as $(1+x_{1}d_{1}+x_{2}a)(1+x_{1}d_{2}+x_{2}b)(1+x_{1}d_{3}+x_{2}c)-f(x_{1},x_{2})$.
Therefore we have the following set of equations in $d,e,f$:
\begin{align*}
&d^{2}+e^{2}+f^{2}=3.25 \\
&d^{2}+2e^{2}+3f^{2}=4.5 \\
&2.5d^{2}+4e^{2}+4.5f^{2}-2def=9
\end{align*}
The set $\mathcal{K}$ of solutions of the first two linear equations in $d^{2},e^{2}$ and $f^{2}$ is defined as follows.
\begin{equation*}
\mathcal{K}:=\{\bmatrix{d^{2}\\e^{2}\\f^{2}}\}=\{\bmatrix{2.5\\.25\\.5}+k\bmatrix{1\\-2\\1}\}
\end{equation*}
Using the third equation we derive the following cubic equation.
\begin{equation*}
8 k^{3}+24 k^{2}+6 k-1=0
\end{equation*}
It provides us $k=-2.7057,-.4076,.1133$.
Solving these polynomial equations we can derive $8$ monic symmetric representations and two of them are non-equivalent. So, two non-equivalent representations of the possible solutions are
\begin{enumerate}
\item $d=-1.616658,e=.152704,f=-.783161$ and $A_{2}=\bmatrix{4.5 & -1.6166 & .1527\\-1.6166 & 4 & -.7831\\.1527 & -.7831 & 2.5}$ at $k=.1133$.
\item $d=-1.446512,e=-1.032089,f=.303968$ and $A_{2}=\bmatrix{4.5 & -1.4465 & -1.0321 \\-1.4465 & 4 & .3040 \\ -1.0321 & .3040 & 2.5}$ at $k=-.4076$.
\item $k=-2.7057$ cannot provide solution as $d^{2}$ is negative in this case which is not possible.
\end{enumerate}
\normalsize
}
\end{example}
\begin{remark} \rm{
While doing the calculation in SAGE, we notice that there are solutions which produce equivalent representations and rest provide non-equivalent representations. During our experiment we derive exactly $2^2.3=12$ solutions. If MSDR of size $3$ exists for a cubic bivariate polynomial, there are $2^2.2=8$ real definite representations and $2$ of them are non-equivalent real definite out of $12$ solutions when the diagonal matrices are simple, the first case. 
}
\end{remark}

\subsection*{Quartic Bivariate Polynomials}
We consider the case of quartic bivariate polynomial and explain the method symbolically as well as numerically. Consider a determinantal  quartic bivariate polynomial
 \begin{align*}
 f(\x):=&f_{40}x_{1}^{4}+f_{31}x_{1}^{3}x_{2}+f_{22}x_{1}^{2}x_{2}^{2}+f_{13}x_{1}x_{2}^{3}+f_{30}x_{1}^{3}+
 f_{21}x_{1}^{2}x_{2}+f_{12}x_{1}x_{2}^{2}\\&+f_{20}x_{1}^{2}+f_{11}x_{1}x_{2}+f_{02}x_{2}^{2}+f_{10}x_{1}+f_{01}x_{2}+1
 \end{align*}
%

For example, by the Theorem \ref{themgmd} we have the following relations between the coefficients of $f(\x)$ and the entries of $D_{1}$ and $A_{12}$ for quartic bivariate case.
 \begin{align*}
 f_{02}&=ab+ac+ad+bc+bd+cd-e^2-f^2-g^2-h^2-k^2-l^2 \\
 f_{12}&=d_{1}(bd+cd+bc)+d_{2}(ad+cd+ac)+d_{3}(ad+bd+ab)+d_{4}(ab+ac+bc)
 -\\&e^2(d_{3}+d_{4})-f^2(d_{2}+d_{4})-g^2(d_{2}+d_{3})-h^2(d_{1}+d_{4})-k^2(d_{1}+d_{3})-l^2(d_{1}+d_{2})\\ f_{22}&=d_{1}d_{2}(cd-l^2)+d_{1}d_{3}(bd-k^2)+d_{1}d_{4}(bc-h^2)+d_{3}d_{4}(ab-e^2)+d_{2}d_{3}(ad-g^2)\\&+d_{2}d_{4}(ac-f^2) \\
 f_{03}&=(abc+acd+abd+bcd)-e^2(c+d)-f^2(b+d)-g^2(b+c)-h^2(a+d)\\&-k^2(a+c)
 -l^2(a+b)+2(hkl+fgl+egk+efh) \\
 f_{13}&=d_{1}(bcd+2hkl)+d_{2}(acd+2fgl)+d_{3}(abd+2egk)+d_{4}(abc+2efh)-e^2(d_{3}d+d_{4}c)\\
 &-f^2(d_{2}d+d_{4}b)-g^2(d_{2}c+d_{3}b)-h^2(d_{4}a+d_{1}d)-k^2(d_{3}a+d_{1}c)-l^2(d_{1}b+d_{2}a)\\
 f_{04}&=abcd-cde^2-bdf^2-bcg^2-adh^2-ack^2-abl^2+2(fglb+hkla+egkc+efhd)\\
 &+e^2l^2+f^2k^2+g^2h^2-2(efkl+eghl+fghk)
 \end{align*}
By comparing the expressions of the coefficients of monomials
 $x_{2}^{2},x_{1} x_{2}^{2}, x_{1}^{2}x_{2}^{2},x_{2}^{3}, x_{1}x_{2}^{3},x_{2}^{4}$ we  obtain an ideal $\I$ generated by three quadratic,  two cubic and one quartic equations. 

In order to compute $A_{12}$ one needs to find the off-diagonal entries of $A_{12}$ that can be obtained  by finding an element of the real variety $V_{\R}(\I)$. Now we see this by a numerical example.
\begin{example}\rm{
Consider a Helton-Vinnikov curve
\begin{equation*}
f(x_{1},x_{2})=1/2 x_{1}^{4}+1/2 x_{2}^{4}-1.5 x_{1}^{2}-1.5 x_{2}^{2}+1/2 x_{1}^{2} x_{2}^{2}+1
\end{equation*}
whose homogeneous version is given in the paper \cite{Sturmfelsquartic}.
The diagonal entries of $D_{1}$ and $D_{2}$ are $(1, 1/\sqrt{2},-1/\sqrt{2},-1)$ and $(0,0,0,0)$ respectively. By the Theorem \ref{themgmd} derive analytic expressions of the coefficients of $f(\x)$ in terms of the coefficient matrices 
\begin{equation*}
 D_{1}=\bmatrix{1 & 0 & 0 &0\\0 & 1/\sqrt(2)& 0&0\\0 & 0 & -1/\sqrt(2) &0 \\0 & 0 & 0 & -1},A_{12}=\bmatrix{0 & e & f & g \\ e & 0 & h & k \\ f & h & 0 & l \\ g & k & l & 0}.
 \end{equation*}
 as follows:
\begin{align*}
f_{02}=&-e^{2}-f^{2}-g^{2}-h^{2}-k^{2}-l^{2} \\
f_{12}=&1.7071e^{2}+.2929f^{2}-.2929k^{2}-1.7071l^{2} \\
f_{22}=&-.7071e^{2}+.7071f^{2}+.5g^{2}+h^{2}+.7071k^{2}-.7071l^{2} \\
f_{03}=&2(hkl+fgl+egk+efh) \\
f_{13}=&2hkl+1.4142fgl-1.4142egk-2efh \\
f_{04}=&e^{2}l^{2}+f^{2}k^{2}+g^{2}h^{2}-2(efkl+eghl+fghk)
\end{align*}
Solving the system of six equations in six unknowns we derive one of the possible representations as follows:
\begin{equation*}
f(x_{1},x_{2})=\det(I+x_{1}D_{1}+x_{2}\bmatrix{0 & .4631 & 0 & .7318\\ .4631 & 0 & -.7318 & 0\\0 & -.7318 & 0 & .4631 \\ .7318 & 0 & .4631 & 0}
\end{equation*}
}
\end{example}

A bivariate polynomial $f(\x)$ of degree $d$ admits an MHDR of size $d$ i.e., there exists an unitary matrix $U_{12}$
\begin{equation}
 f(\x)=\det(I+x_{1}D_{1}+x_{2}U_{12}D_{2}U_{12}^{\ast})
\end{equation}
As it is known that for a bivariate polynomial, MHDR exists if and only if MSDR exists \cite{Helton}, so we focus on constructing MSDR first. One can compute MHDR of size $d$ by solving a different set of polynomial equations, although non-emptyness of $\V_{\R}(\I)$ is a necessary and sufficient condition for the existence of MHDR of size $d$.
\subsection{Equivalent and Non-equivalent Representations} \label{equivclass}
It is evident that if a polynomial $f(\x)$ is determinantal, there are infinitely many such representations defined in equation (\ref{eqdet}). So, it would be an interesting problem to know how many of them are orthogonally or unitarily equivalent to each other.

Observe that two MSDRs of a determinantal bivariate polynomial are orthogonally equivalent i.e., $f(\x)=\det(I+x_{1}D_{1}+x_{2}V_{1}D_{2}V_{1}^{T})=\det(I+x_{1}D_{1}+x_{2}V_{2}D_{2}V_{2}^{T})$ if and only if either $V_{1}D_{2}V_{1}^{T}=V_{2}D_{2}V_{2}^{T}$ or the diagonal matrix $D_{1}$ is invariant under pre and post multiplication by the orthogonal matrix $V_{1}V_{2}^{T}$.

Let matrix $D_{\pm}$ denotes a diagonal matrix with diagonal entries being $\pm 1$, known as \textit{signature matrix}.
\begin{lemma} \label{lemmadiaginv}
If all the eigenvalues of a diagonal matrix $D$ are distinct, then matrix $D$ is invariant under pre and post multiplication by an orthogonal matrix $W$ if and only if ~ $W=D_{\pm}$ is a \textit{signature matrix}.
\end{lemma}
\pf   A diagonal matrix $D$  is invariant under pre and post multiplication by the orthogonal matrix, i.e., $W D W^{T}=D$   if and only if $D$ and $W$ commute  if and only if they have a common eigenvectors. Since all the eigenvalues of $D$ are distinct, so  the corresponding  non-degenerate eigenvectors are unique up to the  sign. Therefore, orthogonal matrix $W$  must be a signature matrix, i.e.,  $W=D_{\pm}$.  \qed

Therefore, the equivalence class of an orthogonal matrix $V_{12}$ defined in equation (\ref{eqbivdet}) is given by $\{D_{\pm}V_{12}\}$. Thus the number of non equivalent MSDR depends on the fact that how many different coefficient matrices $A_{12}$ can be built up using different orthogonal matrices $V_{12}$ up to equivalence class of $\{D_{\pm}V_{12}\}$. 

Note that the diagonal matrices remain unchanged and only the off diagonal entries of a matrix are changed by sign under the action of pre and post multiplication by a signature matrix. So, there are  exactly $2^{d-1}$ coefficient matrices $A_{12}$ with same $D_{1}$ in an equivalent class of MSDRs.

Although the cardinality of the set of equivalence class $\{D_{\pm}V_{12}\}$ is $2^{d}$, but the signature matrix with all diagonal entries $-1$ provides the same coefficient matrix $A_{12}$ like the signature matrix with all diagonal entries $1$ gives.

\begin{remark}\label{remarknoneqvcubic}\rm{
The number of orthogonally non equivalent MSDRs is equal to the number of different (i.e., not equivalent by signature matrix) coefficient matrices $A_{12}$ 
in equation (\ref{eqbivdet}).
}
\end{remark}

It is proved that for smooth Helton-Vinnikov plane curve the number of real definite equivalence classes equals $2^{d-1 \choose 2}$ \cite{Vinnikovselfadjoint},\cite{Sturmfelsbivariate}. 

Using this result we can analyze that for a smooth cubic bivariate polynomial there are $2$ non-equivalent representations and each of these $2$ non-equivalent representation has $2^2=4$ representatives which produce equivalent MSDRs with same diagonal matrix $D_{1}$.

In the same line of thoughts, we conclude that
\begin{corollary}
 If all the eigenvalues of a diagonal matrix $D$ are distinct, then matrix $D$ is invariant under pre and post multiplication by a unitary matrix $W$ if and only if ~ $W=U_{\pm}$ where
 \begin{equation*}
 U_{\pm} =\bmatrix{e^{i \psi_{1}} & 0 & 0\\0 & e^{i \psi_{2}}\\ \vdots & \dots & \vdots \\0 & 0 & e^{i \psi_{d}}} 
\end{equation*}
is a generalization of signature matrix $D_{\pm}$, and we call it \textit{complex signature matrix}.
\end{corollary}
Therefore, the equivalence class of a unitary matrix $U_{12}$ defined in equation (\ref{eqbivdet}) is given by $\{U_{\pm}U_{12}\}$. 

Under the action of pre and post multiplication by a unitary matrix $U_{\pm}$ as the diagonal matrices remain unchanged and only the off diagonal entries of a Hermitian matrix are changed by phases, so it provides infinitely many equivalent coefficient matrices $A_{12}$ with same $D_{1}$ in a equivalent class of MHDRs. 

Thus two unitarily equivalent MHDRs in equation (\ref{eqbivdet}) are same up to the phases of the off diagonal entries of $A_{12}$.

We state the difference between computing a MSDR and MHDR by examples for cubic and quartic bivariate polynomials. The method can be generalized by choosing suitable phases for bivariate polynomials of degree $d$.

For a cubic bivariate polynomial by choosing the phases $(\theta_{1},\theta_{2},\theta_{3})$ such that $\cos(\theta_{1}+\theta_{3}-\theta_{2})=1$
we obtain infinitely many MHDRs such that the magnitudes of the off diagonal entries of Hermitian coefficient matrices associated with $x_{2}$ variable is same as the off diagonal entries of the corresponding symmetric coefficient matrix of the given MSDR. In fact, each such an equivalence class of MHDRs contains the corresponding equivalence class of MSDRs.

For a quartic bivariate polynomial by choosing the cosine of the following phases  equal to one;  i.e.,  by choosing $(\ast)$
\begin{align*}
&\cos(\theta_{1}+\theta_{4}-\theta_{2})=1, \cos(\theta_{1}+\theta_{5}-\theta_{3})=1, \cos(\theta_{2}+\theta_{6}-\theta_{3})=1, \cos(\theta_{4}+\theta_{6}-\theta_{5})=1\\
&\cos(\theta_{1}+\theta_{4}+\theta_{6}-\theta_{3})=1, \cos(\theta_{1}+\theta_{5}-\theta_{2}-\theta_{6})=1, \cos(\theta_{2}+\theta_{5}-\theta_{3}-\theta_{4})=1
\end{align*}
we obtain infinitely many MHDRs such that the magnitudes of the off diagonal entries of Hermitian coefficient matrices associated with $x_{2}$ variable  is same as the off diagonal entries of the corresponding symmetric coefficient matrix of the given MSDR.   

Therefore each choice of different phases such that cosine values of the mentioned phases equal to one is associated with one equivalence class of MHDRs and each such equivalence class of an MHDRs contains the corresponding equivalence class of MSDRs. 

Experiments for both of these cases invoke that  there are infinitely many equivalence classes of MHDRs. For example, by choosing $\cos(\theta_{1}+\theta_{3}-\theta_{2}) \in (-1,1)$ and the cosine values of the mentioned phases $(\ast)$ lying between $(-1,1)$ we can generate a different ideal  and  each point of the real variety of this ideal  provides a non-equivalent MHDR of size $3$ and $4$ respectively. Therefore, we conjecture that  

\begin{conjecture}
There is a continuum of unitarily non-equivalent MHDRs for a bivariate polynomial of degree $d$.
\end{conjecture}
\begin{remark}\rm{
 In order to deal with issue of ideal $\I$ being not a zero dimensional ideal and find one representative from each equivalence class another method to compute MSDR of size $d$ of a bivariate polynomial has been introduced in \cite{papribiv}.
 } 
\end{remark}

\section{Multivariate Polynomials}
In this section, we propose a heuristic method to determine an MSDR of size $d$ for a multivariate polynomial of degree $d$. Here we assume the eigenvalues of all the coefficient matrices are distinct which is equivalent to the fact that the plane curve defined by bivariate polynomial is smooth. The necessity of this assumption is explained later. 

In order to determine an MSDR of size $d$ for a multivariate polynomial of degree $d$ we do the following.
\begin{itemize}
\item Generate $n$ univariate polynomial $f_{x_{i}}:=f(0,\dots,x_{i},\dots,o)$ from the given polynomial by taking its restriction along each of $x_{i},i=1,\dots,n$-th coordinates.
By Lemma \ref{existence of diagonal 2}  if a multivariate polynomial $f(\x) \in \R[\x]$ of degree $d$ has an MSDR (MHDR) of size $d$, then all the roots of $f_{x_{i}}$ are real for all $i=1,\dots,n$. This provides a necessary condition for the existence of an MSDR of a multivariate polynomial.
\item Determine the eigenvalues of coefficient matrices by using Lemma \ref{lemmartdiag}.
The entries of diagonal matrices $D_{i}$ are the negative reciprocal of the roots of univariate polynomials $f_{x_{i}}$ for all $i=1, \dots,n$.  The entries of $D_{i}$ are the eigenvalues of coefficient matrices $A_{1i}$ for all $i=2,\dots,n$. Say $D_{1}=\diag(r_{1},r_{2},\dots,r_{d})$ .
\item  From a given $n$-variate polynomial construct $n \choose 2$ bivariate polynomials $f_{x_{i}x_{j}},i,j=1,\dots,n,i<j$ by making $n-2$ variables zero at a time. 
 Determine an MSDR of size $d$ for those $n \choose 2$ bivariate polynomials by the Algorithm $1$ and then check the compatibility condition which is discussed later.
\item Find a suitable tuple $(D_{1},A_{12},A_{13},\dots,A_{1n})$ using compatibilty condition. Then one needs to compute the determinant of the LMP associated with the suitable tuple 
$(D_{1},A_{12},A_{13},\dots,A_{1n})$. If it coincides with the given polynomial, that tuple provides an MSDR of size $d$.  

Note that determinantal representation for a bivariate polynomial may not be unique up to equivalence class. Therefore, suitable tuple may not be unique for higher degree multivariate polynomials. Though it is observed to be unique for cubic multivariate (more than $2$ variables) case. When the number of variables increases, the number of suitable tuples decreases even though the degree of the polynomial is higher than $3$.

We need to check the determinant of all of these LMPs and if none of them coincides with the given polynomial, conclude that multivariate polynomial of degree $d$ does not admit an MSDR of size $d$. 
\end{itemize}
\subsection*{Compatibility Conditions:}
We need to define compatibility conditions among transition matrices to choose a suitable $n$-tuple of $(V_{12}, V_{13}, \dots, V_{1n})$. Using the method to determine MSDR for a bivariate polynomial we find the orthogonal matrices $V_{1i},V_{1j},V_{ij}$ such that the polynomials  $f_{x_{1}x_{i}}=\det(I+x_{1}D_{1}+x_{i}V_{1i}D_{i}V_{1i}^{T}), f_{x_{1}x_{j}}=\det(I+x_{1}D_{1}+x_{j}V_{1j}D_{j}V_{1j}^{T}),$ and $f_{x_{i}x_{j}}=\det(I+x_{i}D_{i}+x_{j}V_{ij}D_{j}V_{ij}^{T})$. Observe that $\det(I+x_{i}D_{i}+x_{j}V_{ij}D_{j}V_{ij}^{T})=\det(I+x_{i}V_{1i}D_{i}V_{1i}^{T}+x_{j} V_{1i}V_{ij}D_{j}V_{ij}^{T}V_{1i}^{T})$. So, by Lemma \ref{lemmadiaginv} the triplet $(V_{1i},V_{1j}, V_{ij})$ is compatible if and only if $V_{1j}=D_{\pm}V_{1i}V_{ij}$.

However, it is shown  in \cite{plaumann2} that  if the eigenvalues of all coefficient matrices are distinct,  finitely many orthogonally non-equivalent MSDRs of size $d$  are possible for a bivariate polynomial of degree $d$ and the number of non-equivalent real definite representations is precisely $2 ^{d-1 \choose 2}$. 

On the other hand, if a multivariate polynomial $f(\x)$ is determinantal, $n \choose 2$ bivariate polynomials $f_{x_{i}x_{j}},i,j=1,\dots,n,i<j$ are also determinantal, though the converse may not true. So, it is evident that there are finitely many orthogonally non-equivalent MSDRs for a multivariate polynomial if the eigenvalues of all the coefficient matrices are distinct. So, the assumption of distinct eigenvalues of all coefficient matrices make sure that the process of testing compatibility condition will end at finite steps.  

The compatibility condition on transition matrices  is not sufficient condition for the existence of an MSDR for a multivariate polynomial as it does not ensure to satisfy the coefficients of monomials in more than two variables.  
In order to find  compatible $n$-tuple  $(D_{1},A_{12}, \dots,A_{1n})$ of coefficient matrices of a determinantal representation of a multivariate polynomial we use the following iterative process. 

In order to determine a determinantal representation for the part of $f(\x)$ which is a $d$ degree polynomial in three variables  we want to find a $3$ tuple compatible coefficient matrices. Without loss of generality we choose one possible  coefficient matrix $A_{12}$ out of $2^{d-1 \choose 2}$ choices  and do the experiment as follows. 

Suppose we want to  find $A_{13}$ such that $(D_{1},A_{12},A_{13})$ provides a determinantal representation for polynomial $f_{x_{1}x_{2}x_{3}}:=f(x_{1},x_{2},x_{3},0,\dots,0)$ of $f(\x)$ . Note that the diagonal entries of matrix $A_{13}$  are uniquely determined  by  Lemma \ref{lemmauniquediag}. We determine the set of $d \choose 2$ off-diagonal entries of $A_{13}$ by solving another system of linear equations associated with the coefficients of  monomials $x_{2}x_{3},\dots,x_{2}^{d-1}x_{3},x_{1}^{\alpha_{1}}x_{2}^{\alpha_{2}}x_{3}$, $1 \leq \alpha_{1} \leq d-2,1 \leq \alpha_{2} \leq d-2,\alpha_{1}+\alpha_{2} \leq d-1$. The number of such monomials is $d \choose 2$. If the system of $d \choose 2$ linear equations in $d \choose 2$ variables  is inconsistent, start with another $A_{12}$. 

Otherwise, we repeat the same method by including one more variable at each step to get a $n$-tuple compatible coefficient matrices. Note that generically the set of off-diagonal entries of $A_{13}$ is uniquely determined when the eigenvalues of coefficient matrices are distinct for a fixed $A_{12}$. If the method fails for all possible coefficient matrix $A_{12}$, we declare that the multivariate polynomial $f(\x)$ has no MSDR of size $d$.

\subsection{Cubic Multivariate Determinantal Polynomials} \label{secmsdrcubic}

First we deal with this issue for cubic trivariate polynomials. In order to determine an MSDR of a trivariate polynomial we need to find orthogonal matrices $V_{12}, V_{13}$ such that $f(\x)=\det(I+x_{1}D_{1}+x_{2}V_{12}D_{2}V_{12}^{T}+x_{3} V_{13}D_{3}V_{13}^{T})$.

Consider the cubic trivariate polynomial $f(\x) \in \R[x_{1},x_{2},x_{3}]$ such that
\begin{eqnarray*}
f(\x)&=f_{300}x_{1}^{3}+f_{210}x_{1}^{2}x_{2}+f_{120}x_{1}x_{2}^{2}+f_{201}x_{1}^{2}x_{3}+f_{102}x_{1}x_{3}^{2}+f_{021}x_{2}^{2}x_{3}\\&+f_{111}x_{1}x_{2}x_{3}+
f_{012}x_{2}x_{3}^{2}+f_{030}x_{2}^{3}+f_{003}x_{3}^{3}+f_{200}x_{1}^{2}+f_{110}x_{1}x_{2}+f_{101}x_{1}x_{3}\\&+f_{011}x_{2}x_{3}+f_{020}x_{2}^{2}+f_{002}x_{3}^{2}
+f_{100}x_{1}+f_{010}x_{2}+f_{001}x_{3}+1.
\end{eqnarray*}
Indeed the vector coefficient of the mixed monomials are dependent on the matrices $D_{1}, A_{12}:=V_{12}D_{2}V_{12}^{T},$ and $A_{13}:=V_{13}D_{3}V_{13}^{T}$. Observe that if $f(\x)=\det(I+x_{1}D_{1}+V_{12}D_{2}V_{12}^{T}+x_{3} V_{13}D_{3}V_{13}^{T})$, then bivariate polynomials $f_{x_{i},x_{j}}$ admit MSDRs of the form $f_{x_{i},x_{j}}=\det(I+x_{i}D_{i}+V_{1j}D_{j}V_{1j}^{T})$ for all $i,j=1,2,3,i <j$. 

Note that any three tuple $(V_{12}, V_{13}, V_{23})$ of orthogonal matrices such that $f(x_{1},x_{2})=\det(I+x_{1}D_{1}+x_{2}V_{12}D_{2}V_{12}^{T}), f(x_{1},x_{3})=\det(I+x_{1}D_{1}+x_{3}V_{13}D_{3}V_{13}^{T})$ and $f(x_{2},x_{3})=\det(I+x_{2}D_{2}+x_{3}V_{23}D_{3}V_{23}^{T})$ does not ensure that $f(x_{1},x_{2},x_{3})=\det(I+x_{1}D_{1}+x_{2}V_{12}D_{2}V_{12}^{T}+x_{3}V_{13}D_{3}V_{13}^{T})$. So, the converse statement need not be true. The converse statement will be true if we can find a suitable combination of $V_{12}, V_{13}$ such that $A_{12}:=V_{12}D_{2}V_{12}^{T} \in \mathcal{O}_{D_{2}}$, and $A_{13}:=V_{13}D_{3}V_{13}^{T} \in \mathcal{O}_{D_{3}}$ would be the required coefficient matrices.

Observe that $f_{x_{2},x_{3}}=\det(I+x_{2}D_{2}+x_{3} V_{23}D_{3}V_{23}^{T})$ and $\det(I+x_{2}D_{2}+x_{3}V_{23}D_{3}V_{23}^{T})=\det(I+x_{2}V_{12}D_{2}V_{12}^{T}+x_{3} V_{12}V_{23}D_{3}V_{23}^{T}V_{12}^{T})$. So, if the required $3$ tuple $(V_{12},V_{13}, V_{23})$ of orthogonal matrices is compatible, then $V_{13}=V_{12}V_{23}$ in other words, the following diagram commutes.

\begin{displaymath}
    \xymatrix{
       \mathcal{O}_{D_{1}} \ar[d]^{V_{13}}  \ar[r]^{V_{12}} & \mathcal{O}_{D_{2}} \ar[dl]^{V_{23}} \\
        \mathcal{O}_{D_{3}}                      &  }
\end{displaymath}

The following theorem provides a necessary and sufficient condition for the existence of MSDR of size $3$ for cubic trivariate polynomials.
\begin{theorem} \label{compatibilitythm}
A cubic trivariate polynomial $f(\x)$ has an MSDR of size $3$ given by $f(\x)=\det(I+\x_{1}D_{1}+\x_{2}V_{12}D_{2}V_{12}^{T}+\x_{3}V_{13}D_{3}V_{13}^{T})$, if and only if the following condition hold:
Existence of MSDR of $3 \choose 2$ bivariate polynomials $f_{(x_{i},x_{j})},i,j=1,2,3,i<j$ (which are constructed from the given polynomial by making one variable zero at a time), and the coefficient matrices of such MSDR of $3 \choose 2$ bivariate polynomials should satisfy the compatibility condition given as follows:
\begin{align*}
A_{13}V_{12}=V_{12}A_{23} \ \mbox{where} \ A_{12}=V_{12}D_{2}V_{12}^{T}, A_{13}=V_{13}D_{3}V_{13}^{T}, \ \mbox{and} \ A_{23}=V_{23}D_{3}V_{23}^{T}.
\end{align*}
for a triplet of orthogonal matrices $(V_{12}, V_{13}, V_{23})$ where $V_{ij}, i \neq j$ is the transition matrix from $D_{i} \in \mathcal{O}_{D_{i}}$ to $A_{ij}:= V_{ij} D_{j}V_{ij}^{T} \in \mathcal{O}_{D_{j}}, i=1,2;j=2,3$ and the pair of coefficient matrices $(A_{12},A_{13})$ satisfy the coefficient of monomial $x_{1}x_{2}x_{3}$ whose analytic expression is given by
\begin{equation*}
f_{111}:=d_{1}(bn+mc-2fq)+d_{2}(an+cl-2ep)+d_{3}(am+bl-2od).
\end{equation*}
\end{theorem}
\pf It is evident that if MSDR of size $3$ exists for a cubic trivariate polynomial, these conditions are satisfied. Conversely, it is clear from the above discussion that if there exists a triplet of orthogonal matrices $(V_{12},V_{13},V_{23})$ such that $f(x_{1},x_{2})=\det(I+x_{1}D_{1}+x_{2}V_{12}D_{2}V_{12}^{T}), f(x_{1},x_{3})=\det(I+x_{1}D_{1}+x_{3}V_{13}D_{3}V_{13}^{T}), f(x_{2},x_{3})=\det(I+x_{2}D_{2}+x_{3}V_{23}D_{3}V_{23}^{T})$ and $V_{13}=V_{12}V_{23}$, then the pair of coefficient matrices $(A_{12},A_{13})$ is compatible. On the other hand, $V_{13}=V_{12}V_{23} \Leftrightarrow A_{13}=V_{12}V_{23}D_{3}V_{23}^{T}V_{12}^{T}=V_{12}A_{23}V_{12}^{T} \Leftrightarrow A_{13} V_{12}=V_{12}A_{23} (\Leftrightarrow D_{3} V_{13}^{T}V_{12}V_{23}=V_{13}^{T}V_{12}V_{23}D_{3})$. Although this is a compatibility condition, but it is not sufficient condition for the existence of an MSDR for a cubic trivariate polynomial as it does not ensure to satisfy the coefficient of monomial $x_{1}x_{2}x_{3}$ for a cubic trivariate polynomial. So, if this compatible pair $(A_{12},A_{13})$ of coefficient matrices satisfy the coefficient of $x_{1}x_{2}x_{3}$ of given polynomial, then it ensures that the cubic trivariate polynomial $f(\x)$ admits an MSDR of size $3$.  \qed

\begin{remark}\rm{
Observe that the steps in the iterative process of this method is finite as the number of equivalent class is finite \cite{Sturmfelsbivariate}. So,  if MSDR for cubic trivariate polynomial exists, this method must provide the required coefficient matrix $A_{13}$ for a fixed $A_{12}$, otherwise MSDR of size $3$ does not exist.
}
\end{remark}
\begin{remark}\rm{
In fact, there are  $2$ non equivalent MSDRs for a smooth cubic plane curve defined by bivariate polynomial, therefore one has to verify $4$ cases to determine MSDRs. Fix $A_{12}$, find a compatible $A_{13}$. If $A_{13}$ exists for one chosen coefficient matrix $A_{12}$, then it happens for all such coefficient matrices (associated to $x_{2}$ variable) which are equivalent to $A_{12}$. 
}
\end{remark}
Observe that for any two non equivalent MSDRs, the diagonal entries of the coefficient matrices $A_{12}$ and $A_{13}$ are invariant as they are obtained by solving a system of linear equations after fixing the diagonal matrix $D_{1}$. Using this fact we propose one more way to determine the off-diagonal entries of the coefficient matrix $A_{13}$ after fixing the coefficient matrix $A_{12}$ to avoid to check compatibility condition directly.
Suppose
\begin{align*}
f(x_{1},x_{2},x_{3})&=\det(I+x_{1}D_{1}+x_{2}A_{12}+x_{3}A_{13})\\
&\det(I+x_{1}\bmatrix{d_{1} & 0 &0\\0 &d_{2}&0\\0&0&d_{3}}+x_{2}\bmatrix{a & d& e\\d&b&f\\e&f&c}+x_{3}\bmatrix{l & o &p\\o &m&q\\p&q&n})
\end{align*}
After evaluating the values of coefficient matrices $D_{1},A_{12}$ and the diagonal entries of $A_{13}$, choose one of the possible coefficient matrices $A_{12}$. Then observe that the coefficients of monomials $x_{2}x_{3},x_{2}^{2}x_{3},x_{1}x_{2}x_{3}$ can be expressed as linear equations in terms of off diagonal entries of $A_{13}$ as follows.
\begin{align*}
f_{011}&=a(m+n)+b(l+n)+c(l+m)-2(od+pe+fq)\\
f_{021}&=l(bc-f^2)+m(ac-e^2)+n(ab-d^2)+2p(df-be)+2q(de-af)+2o(ef-cd)\\
f_{111}&=d_{1}(bn+mc-2fq)+d_{2}(an+cl-2ep)+d_{3}(am+bl-2od)
\end{align*}
By solving a system of linear equations of the form $H \y=\z$, we can determine the off diagonal entries $o,p,q$ of coefficient matrix $A_{13}$, where
\tiny
\begin{equation} \label{eqcubictriv}
H=2\bmatrix{d & e & f\\(cd-ef) & (be-df) & (af-de)\\d_{3}d & d_{2}e & d_{1}f}, \y=\bmatrix{o\\p\\q}, \z=\bmatrix{a(m+n)+b(l+n)+c(l+m)-f_{011}\\l(bc-f^{2})+m(ac-e^{2})+n(ab-d^{2})-f_{021}\\d_{1}(bn+mc)+d_{2}(an+cl)+d_{3}(am+bl)-f_{111}}.
\end{equation}
\normalsize
If $\z \in \mbox{col}(H)$, then the solution of this system of linear equations exists and if rank($H)=3$, then unique solution exists.
Next check whether it satisfies the expression due to the coefficients of monomials $x_{1}x_{3}^{2}, x_{2}x_{3}^{2}, x_{3}^{2},x_{3}^{3}$ as follows.
\begin{align} \label{cubictriveq2}
f_{102}&=d_{1}mn+d_{2}ln+d_{3}lm-d_{1}q^{2}-d_{2}p^{2}-d_{3}o^{2} \nonumber \\
f_{012}&=a(mn-q^{2})+b(ln-p^{2})ln+c(lm-o^{2})+2d(pq-on)+2e(oq-mp)+2f(op-lq) \nonumber \\
f_{002}&=mn+ln+lm-o^{2}-p^{2}-q^{2} \nonumber \\
f_{003}&=lmn+2opq-lq^{2}-mp^{2}-no^{2} \nonumber 
\end{align}
This can be checked by constructing matrix $A_{13}$ and calculating the determinant of LMP. If the system has infinitely many solutions, the set of solutions
can be written as parameterized solutions in one or two parameters depending on the rank of matrix $H$. Then eliminate those parameters using the equations associated with $f_{102},f_{012},f_{002},f_{003}$ mentioned above and determine the off diagonal entries of $A_{13}$. 

We see an example.
\begin{example}\rm{
Consider the cubic trivariate polynomial 
\begin{eqnarray*}
&1+6x_{1}^3+36x_{1}^2x_{2} + 66x_{1}x_{2}^2 + 36x_{2}^{3}+70x_{1}^{2}x_{3}+210x_{1}x_{3}^{2}+ 162x_{3}^{3}+366.819x_{2}x_{3}^{2}+225.7077x_{2}^{2}x_{3}\\&+262.2732 x_{1}x_{2}x_{3}+ 11x_{1}^{2}+42x_{1}x_{2}+36x_{2}^2 +74x_{1}x_{3}+99x_{3}^{2}+133.1368x_{2}x_{3}+6 x_{1}+ 11x_{2} + 18 x_{3}
\end{eqnarray*}
By Lemma \ref{lemmartdiag} $D_{1}=\diag(3,2,1),D_{2}=\diag(6,3,2)$ and $D_{3}=\diag(9,6,3)$. Two non-equivalent MSDRs (shown in Example \ref{excubicbivtensor}) for bivariate polynomial
 \begin{equation*}
  f_{(x_{1},x_{2})}=6x_{1}^3 + 36x_{1}^2x_{2} + 11x_{1}^2 + 66x_{1}x_{2}^2 + 42x_{1}x_{2} + 6x_{1} + 36x_{2}^3 + 36x_{2}^2 + 11x_{2} + 1
  \end{equation*}
  are gievn by (one possible representative from each equivalence class)
 \begin{equation*}
  A_{12}=\bmatrix{4.5 & -1.6166 & .1527\\-1.6166 & 4 & -.7831\\.1527 & -.7831 & 2.5},
 A_{12}=\bmatrix{4.5 & -1.4465 & -1.0321 \\-1.4465 & 4 & .3040 \\ -1.0321 & .3040 & 2.5} 
 \end{equation*}
For bivariate polynomial
\begin{equation*}
f_{(x_{1},x_{3})}=6x_{1}^{3}+70x_{1}^{2}x_{3}+210x_{1}x_{3}^{2}+74x_{1}x_{3}+11x_{1}^{2}+6x_{1}+162x_{3}^{3}+99x_{3}^{2}+18x_{3}+1
\end{equation*}
two non-equivalent MSDRs are given by
\begin{equation*}
A_{13}=\bmatrix{ 5 & 0 & 2.8284\\0 & 6 & 0 \\ 2.8284 &0 &7}, A_{13}= \bmatrix{5 & -2 & 0\\-2 & 6 & 2\\0 & 2 & 7}.
\end{equation*}
We have four possibilities to find a suitable tuple $(D_{1},A_{12},A_{13})$ such that the determinant of corresponding LMP coincides with the given trivariate polynomial. Instead of trying all possible situations (which is $2^{d}$ for degree $d$ polynomial) we make use of the compatibility condition. We find that the (unique) compatible tuple is 
\begin{equation*}
(D_{1}=\diag(3,2,1),A_{12}=\bmatrix{4.5 & -1.6166 & 0.1527\\ -1.6166 & 4 & -0.7831\\ 0.1527 & -0.7831 & 2.5}, A_{13}=\bmatrix{5& 0& 2.8284\\0 & 6 & 0\\ 2.8284 & 0 &7})
\end{equation*}
In the alternative method, we compute $D_{1},D_{2},D_{3},A_{12}$ and $\diag(A_{13})$. Then we compute the off-diagonal entries of suitable $A_{13}$ by solving a system of linear equations of the form $H \y=\z$ where 
\begin{equation*}
H=\bmatrix{-3.2332 & 0.3054 &-1.5662\\-7.8438 &-1.3103&-6.5542\\-3.2332&   0.6108&-4.6986}, \z=\bmatrix{0.8632\\-3.7076\\1.7268} \ \mbox{by equation (\ref{eqcubictriv})}.
\end{equation*}
As $H$ is nonsingular, so $H \y=\z$ has unique solution i.e., the off diagonal elements of $A_{13}$ are $0,2.8284,0$. Indeed, both of these methods provide a unique equivalent class of MSDR of size $3$ for the given trivariate polynomial.
}
\end{example}
\subsection{Generalization:A Heuristic Method}
Consider an $n$ variate cubic polynomial $f(\x)=\sum_{j_{k} \in \{0,\dots, 3\},\sum_{k=1}^{n} j_{k}\leq 3} f_{j_{1} \dots j_{n}} x_{1}^{j_{1}} \dots x_{n}^{j_{n}}$. Suppose polynomial $f(\x)$ admits an MSDR of size $3$ such that
\begin{align*}
&f(\x)=\det(I+x_{1}D_{1}+x_{2}A_{12}+x_{3}A_{13}+x_{4}A_{14}+\dots+x_{n}A_{1n})\\
&\det(I+x_{1}\bmatrix{d_{1} & 0 &0\\0 &d_{2}&0\\0&0&d_{3}}+x_{2}\bmatrix{a & d & e\\d & b & f\\e & f & c}+x_{3}\bmatrix{l & o &p\\o & m & q\\p & q & n}
+x_{4} \bmatrix{r & u &v\\u&s&w\\v&w&t}+\dots)
\end{align*}
We can realize the complexity of this problem when we deal with four variables. We have been able to determine $D_{1},A_{12},A_{13}$ and diagonal entries of
coefficient matrices $A_{1n}, n > 3$. From the discussion of the previous subsection, it is clear that to determine a suitable $A_{14}$ we have to solve a system of seven linear equations in $3$ unknowns-an overdetermined system. There is a high probability of non existence of an MSDR for a cubic quartic polynomial. These seven linear equations are obtained from the monomials as follows.
\begin{align*} \label{cubicnvariateeq1}
f_{0101}&=a(s+t)+b(r+t)+c(r+s)-2(du+ev+fw)  \\
f_{0201}&=r(bc-f^{2})+s(ac-e^{2})+t(ab-d^{2})+2u(ef-cd)+2v(df-be)+2w(de-af)\\
f_{0011}&=l(s+t)+m(r+t)+n(r+s)-2(ou+pv+qw)  \\
f_{0021}&=r(mn-q^{2})+s(ln-p^{2})+t(lm-o^{2})+2u(pq-on)+2v(oq-mp)+2w(op-pq)\\
f_{1101}&= d_{1}(bt+sc-2fw)+d_{2}(at+cr-2ev)+d_{3}(as+br-2du)\\
f_{1011}&=d_{1}(mt+ns-2qw)+d_{2}(lt+nr-2pv)+d_{3}(ls+mr-2ou)  \\
f_{0111}&=a(mt+ns-2qw)+b(lt+nr-2pv)+c(ls+mr-2ou)+2d(pw+vq-ot-nu)\\&+2e(ow+qu-mv-ps)+2f(ov+pu-lw-qr)
\end{align*}
Therefore, we propose a heuristic method for cubic multivariate (more than $3$ variables)  polynomials based on compatibility condtitions.

We explain some properties of transition matrices which will be used later. Say $V_{ij}, i \neq j$ is the transition matrix from $D_{i} \in \mathcal{O}_{D_{i}}$ to $A_{ij}:= V_{ij} D_{j}V_{ij}^{T} \in \mathcal{O}_{D_{j}}$ and similarly $V_{ij}^{T}=V_{ji}, i \neq j$ is the transition matrix from $D_{j} \in \mathcal{O}_{D_{j}}$ to $ A_{ji}:=V_{ij}^{T} D_{i}V_{ij}=V_{ji}D_{i}V_{ji} \in \mathcal{O}_{D_{i}}$.
\begin{proposition} \label{propcompatibility}
$V_{1j}V_{jk}=V_{1k}$ for all $j,k (j \neq k, j, k \in \{2, \dots, n\})$ if and only if $V_{ij}V_{jk}=V_{ik}$ for all $i,j,k (i \neq j \neq k, i,j,k \in \{1, \dots, n\})$ where $V_{ij}$ is the connecting matrix defined above.
\end{proposition}
\pf R.H.S implies L.H.S is obvious. We need to prove that the condition in L.H.S is sufficient for the condition in R.H.S to be hold.
If $n=3$, there is nothing to prove. To avoid the same parity we have taken the ordering $i < j < k$. Say $n \geq 4$. Aim is to prove that if $n-1 \choose 2$ relations hold, remaining $n \choose 3$ $-$ $n-1 \choose 2$ relations are true. Note that
\begin{align*}
V_{ij}V_{jk}&=V_{i1}V_{1j}V_{jk} (\mbox{as it is given} \ V_{1i}V_{ij}=V_{1j}) \\
&=V_{i1}V_{1k}=V_{ik} (\mbox{as it is given} \ V_{1j}V_{jk}=V_{1k}) \forall i,j,k \in \{2, \dots,n\}.
\end{align*} Hence the proof. \hfill${\square}$

Therefore, using the Theorem \ref{compatibilitythm} and Proposition \ref{propcompatibility} we conclude the following theorem.
\begin{theorem} \label{themmulti}
A multivariate polynomial $f(\x)$ of degree $d$ admits an MSDR of size $d$, given by $f(\x)=\det(I+x_{1}D_{1}+x_{2}V_{12}D_{2}V_{12}^{T}+\dots+x_{n}V_{1n}D_{n}V_{1n}^{T})$, if and only if the following conditions hold:
Existence of MSDR of $n \choose 2$ bivariate polynomials $f_{(x_{i},x_{j})}, i,j=1,\dots,n, i<j$ (which are constructed from the given polynomial by making $n-2$ variables zero at a time), and the coefficient matrices of such MSDR of $n \choose 2$ bivariate polynomials should satisfy the compatibility condition given as follows:
\begin{align*}
A_{1j}V_{1k}=V_{1k}A_{kj}, \forall j=3,\dots,n, k=2,\dots,n-1 .
\end{align*}
for $n-1 \choose 2$ triplet of orthogonal matrices $(V_{1 j-1},V_{1j},V_{j-1 j})$ where $V_{ij}, i \neq j$ is the transition matrix from $D_{i} \in \mathcal{O}_{D_{i}}$ to $A_{ij}:= V_{ij} D_{j}V_{ij}^{T} \in \mathcal{O}_{D_{j}}$ and the $n$-tuple coefficient matrices $(D_{1},A_{12},A_{13}, \dots, A_{1n})$ satisfy all the coefficients of mixed monomials in $3$ to $n$ variables.
\end{theorem}

\bibliographystyle{alpha}
\bibliography{refqplmi}
\end{document}